\documentclass[a4paper]{amsart}
\usepackage{amssymb}

\usepackage[dvips]{graphicx}
\usepackage{latexsym}
\usepackage{amsmath}

\newcommand{\hauteur}{6cm}
\newcommand{\largeur}{10cm}

\theoremstyle{plain}

\numberwithin{equation}{section}

\input{tcilatex}

\begin{document}

\title{On the monotonicity criteria of the period function of potential systems}
\author{A. Raouf Chouikha }
\address{ Universit\'e Paris XIII, Institut Galil\'ee, LAGA CNRS UMR
7539, 93430 Villetaneuse, {\textsc{France}} }
\email{chouikha@math.univ-paris13.fr}
\date{}



\begin{abstract} 
The purpose of this paper is to study various monotonicity conditions of the period function $T(c)$  (energy-dependent) for potential systems $\ddot x + g(x)=0$ with a center at the origin $0$.
We had before identified a family of new criteria noted by $(C_n)$ which are sometimes thinner than those previously known ({\it Period function and characterizations of Isochronous potentials}\quad arXiv:1109.4611). This fact will be illustrated by examples.

{\it Key Words and phrases:} \ period function, monotonicity, center, potential systems.\footnote
{2000 Mathematics Subject Classification  \ 34C05, 34C25, 34C35}. 
\end{abstract}
 
\maketitle

\section {Introduction and statement of results}
Consider the potential system 
\begin{equation} \qquad \dot x = y , \qquad \dot y = - g(x)\end{equation} 
where $\dot x = \frac{dx}{dt}, \ddot x = \frac{d^2x}{dt^2}$ and $g(x)$ is analytic on $R$.\\ Let $G(x)$ be the potential of (1) 
$$G(x) = \int_0^x g(\xi) d\xi.$$
The following hypothesis ensures that (1) has a center at the origin $0$.\\ There  exist  $a < 0 < b$  such  that
  $$(\mathcal{H} )\qquad \left\{ 
   G(a) = G(b) = c,\ G(x) < c \ and \ x g(x) > 0 \ for \  all \ a < x < b \ and \ x \neq 0.
\right\}   $$ {\it Moreover, without lose of generality we will assume in the sequel that\\ $g(0) = 0$ and} $g'(0) = 1.$ \  \\

 Moreover, we consider the involution $A$ defined by (see [5] for example)\\  
$ G(A(x)) = G(x) \ and \ A(x) x < 0 $ for all $x \in [a,b]$. This means $A(0)=0$ and when $x \in [0,b]$ then $A(x) \in [a,0]$.\\
 Let $T(c)$ denotes the minimal period of periodic orbits depending on the energy. 
\begin{equation} T(c) = \sqrt{2}\int_a^b \frac{dx}{\sqrt{c-G(x)}}.\end{equation}  The period function is well defined for any $c$ such that $0 < c < \bar c$ and when $\bar c < +\infty$.

We proved the following (Theorem A of [1]) 
\bigskip

{\bf Theorem 1-1} \quad {\it  Let $g(x)$ be an analytic function, $G(x) = \int_0^x g(\xi) d\xi$ be the potential of equation $(1)\quad \ddot x + g(x) = 0.$ and let $A(x)$ be the involution defined above. Suppose hypothesis \ ($\mathcal{H}$)\ holds, let us define the n-polynomial with respect to $G$, $$f_n(G) = f(0)+f'(0)G+ \frac{1}{2}f''(0)G^2+....+\frac{1}{n!}f^{(n)}(0)G^n$$ where for $g''(0)\neq 0,$ $$ \quad  f(0)= -\frac{g''(0)}{3}=-\frac{1}{3} \frac{d^2 g}{dx^2}(0),\quad f'(0)= -\frac{7}{9}g''^3(0)+ \frac{g^{(4)}(0)}{5},$$ $$ f''(0)=-\frac{28550(g''(0))^6-10320(g''(0))^3g^{(4)}(0)+81(g^{(4)}(0))^2+270(g''(0))g^{(6)}(0)}{4050 g''(0)}
,... $$  (When $g^{(2)}(0)=0,\quad f(0)=0,\quad f'(0)=-\frac{1}{5}g^{(4)}(0),\quad  f''(0)=-\frac{1}{42}g^{(6)}(0),...$)\\ Suppose that for a fixed $n \in N$ and for $x \in [0,b]$ one has $$(\mathcal{C}_n)\qquad \frac{d}{dx}[\frac{G}{g^2}(x)] > f_n(G) > \frac{d}{dx}[\frac{G}{g^2}(A(x))]$$ $$ (or\  \frac{d}{dx}[\frac{G}{g^2}(x)] < f_n(G) < \frac{d}{dx}[\frac{G}{g^2}(A(x))] )$$  then $T(c)$ the period function of (1) is increasing (or decreasing) for $0 < c < \bar c$.}\\ 
 
\bigskip 

{\bf Remark 1} \quad We may also define the coefficients of a simple manner as follows :
$$f(0)= lim_{x\rightarrow 0} \frac{d}{dx}[\frac{G}{g^2}(x)] $$
$$f'(0)=lim_{x\rightarrow 0} \frac{1}{G}[\frac{d}{dx}[\frac{G}{g^2}(x)]-f(0)] $$
$$\frac{f'(0)}{2}=lim_{x\rightarrow 0} \frac{1}{G^2}[\frac{d}{dx}[\frac{G}{g^2}(x)]-f(0)-f'(0)G],.... $$

As a first consequence one deduces the following which has already been proved by Chow-Wang (Cor. 2.5, [4]).

\bigskip
{\bf Corollary 1-2}\quad {\it Suppose hypothesis ($\mathcal{H}$) holds and let $g(x)$ be an analytic function for $ 0 < x <  b $ and $G(x) = \int_0^x g(\xi) d\xi$ be the potential of (1) and $g''(0) \neq 0$.\\
 Suppose condition $(\mathcal{C}_0)$ holds. This means $$\frac{d}{dx}[\frac{G}{g^2}(x)] >(or <) f(0) >(or <) \frac{d}{dx}[\frac{G}{g^2}(A(x))]$$ or equivalently 
$$g^2(x) + \frac {g''(0)}{3} g^3(x) - 2 G(x) g'(x) > 0 (or < 0)$$  then $T(c)$ is increasing (or decreasing) for $0 < c < \bar c$.}
\bigskip

By the same way we may deduce 

\bigskip

{\bf Corollary 1-3}\quad {\it Let $g(x)$ be an analytic function  and \\ $G(x) = \int_0^x g(\xi) d\xi$ be the potential of equation (1). Suppose $(\mathcal{C}_1)$ holds. That means for $0<x<b,$ $$\frac{d}{dx}[\frac{G}{g^2}(x)] >(or <) f(0)+f'(0)G >(or <) \frac{d}{dx}[\frac{G}{g^2}(A(x))]$$ or equivalently $$g^2(x) + \frac {g''(0)}{3 } g^3(x) - 2 G(x) g'(x) + (\frac {7g''^3(0)}{9}- \frac {g^{(4)}(0)}{5}) g^3(x) G(x) > 0 (or < 0)$$ then $T(c)$  is increasing (or decreasing) for $0 < c < \bar c$.  }

\bigskip

{\bf Corollary 1-4}\quad {\it Let $g(x)$ be an analytic function  and \\ $G(x) = \int_0^x g(\xi) d\xi$ be the potential of equation (1). Suppose $(\mathcal{C}_2)$ holds. That means for $0<x<b,$ $$\frac{d}{dx}[\frac{G}{g^2}(x)] >(or <) f(0)+f'(0)G+{\frac{1}{2}f''(0)G^2} >(or <) \frac{d}{dx}[\frac{G}{g^2}(A(x))]$$  then $T(c)$ is increasing (or decreasing) for $0 < c < \bar c$.  }\\
\bigskip

{\bf Remark 2}\quad   The monotonicity problem of the period function has been extensively studied. Many criteria have been produced. A lot of them logically are related. For a comparison between these sufficients conditions  we may refer to [2] and [3] and references therein.\\ 
Although, notice that the monotonicity criterium $(\mathcal{C}_0)$ given by Corollary 1-2 appears sometimes to be the best one. Indeed, it is more general than those given by C. Chicone, F. Rothe [6] and R. Schaaf [7].\\  
In [4] we proved the non-optimality of these criteria by giving appropriate examples of potential $G$ for which the energy-period is monotonic, in spite of none of these conditions of monotony is verified. \\
It thus seems to ask if we could to compare these new conditions each other. We are then content to make a few remark about the sign of $f^{(k)}(0)$.\\
More precisely, it is clear if we suppose $f'(0)= -\frac{7}{9}g''^3(0)+ \frac{g^{(4)}(0)}{5} < 0$ and the potential $G$ satisfies $(\mathcal{C}_0)$  then $G$ also satisfies $(\mathcal{C}_1)$ ( impliyng together $T(c)$ is monotonic). That means $(\mathcal{C}_1)$ is better than $(\mathcal{C}_0)$. We will say in the sequel : "$(\mathcal{C}_0)$\  {\it implies} \ $(\mathcal{C}_1)$". By the same way, when $f''(0) < 0$ then condition $(\mathcal{C}_1)$ \  {\it implies} \ $(\mathcal{C}_2)$. When $f^{3}(0) < 0$ then condition $(\mathcal{C}_2)$ \  {\it implies} \ $(\mathcal{C}_3)$.\\ More generally, we may claim when $f^{(k)}(0) < 0$ then $(\mathcal{C}_{k-1})$ \  {\it implies} \ $(\mathcal{C}_k)$ and when $f^{(k)}(0) > 0$ then $(\mathcal{C}_k)$ \  {\it implies} \ $(\mathcal{C}_{k-1})$.\\

We may ask if these implications are strict. 
Below, we will give an exemple of potential for which condition $(\mathcal{C}_2)$ is verified 
but not $(\mathcal{C}_1)$ nor $(\mathcal{C}_0)$.\\ Before to continue consider at first the following

\section{The case of $g''(0)=0$}

 The pioneering work devoted to the study of the period function is undoubtedly the Opial's paper [6]. He interested in behavior and monotonicity of the period function $T(c)$ of equation (1). When $g''(0)=0$ he proved that condition
 $$
  x\frac{d }{dx} (\frac{g(x) }{x}) \ \neq 0\qquad (Op)
$$
implies $T(c)$ is monotonic.\\ 
 We proved in [2] that the Opial condition $(Op)$ of monotonicity for the period function 
 is the better among all known conditions for which $g''(0) = 0$. Indeed, we may prove the following (which is a slightly modified version of Theorem 3 of [2]) 
\bigskip

 {\bf Theorem 2-1}\  \quad {\it Let $g(x)$ be an analytic function  and \\ $G(x) = \int_0^x g(\xi) d\xi$ be the potential of equation (1) satisfying hypothesis $(\mathcal{H})$. Then we have the following implications\\
$$x g''(x) < 0 \quad {\it implies} \quad g^2(x) - 2G(x) g'(x) > 0 \quad {\it implies} \quad x (\frac{d}{dx}\frac{g(x)}{x})) < 0.$$
Moreover, each of these conditions implies that the period function $T(c)$ of (1) is strictly increasing for $0 < c < c_1$.\\
$$x g''(x) > 0 \quad {\it implies} \quad g^2(x) - 2G(x) g'(x) < 0 \quad {\it implies} \quad x (\frac{d}{dx}\frac{g(x)}{x})) > 0.$$
Moreover, each of these conditions implies that the period function $T(c)$ of (1) is strictly decreasing for $0 < c < c_1$.\\
A necessary condition to have any of these conditions is $g''(0) = 0$ .}

\bigskip

Applying Theorem 1-1, Corollary 1-3 and Remark 2 we prove the following

\bigskip 

{\bf Proposition 2-2} \quad {\it Let $g(x)$ be an analytic function  and \\ $G(x) = \int_0^x g(\xi) d\xi$ be the potential of equation (1) satisfying hypothesis $(\mathcal{H})$. Suppose $g''(0)= 0, g^{(4)}(0) < 0$, then 
$$g^2(x) - 2G(x) g'(x) > 0 (or < 0) \quad {\it implies} \quad (\mathcal{C}_1) \quad {\it implies} \quad x (\frac{d}{dx}\frac{g(x)}{x})) < 0(or > 0) $$ Recall that $$ (\mathcal{C}_1) :\ \frac{d}{dx}[\frac{G}{g^2}(x)] >(or <) \frac{{g^{(4)}}(0)}{5}G >(or <) \frac{d}{dx}[\frac{G}{g^2}(A(x))].$$

Moreover, each of these two conditions implies that the period function $T(c)$ of (1) is strictly increasing (or decreasing) }

\bigskip

{\bf Proof}\quad Indeed, since  $f'(0)= \frac{{g^{(4)}}(0)}{5} < 0$ then by Corollary 1-2 and Remark 2 \ $(\mathcal{C}_0)$ \ {\it implies}\ $(\mathcal{C}_1)$. On the other hand, $$g^2(x) - 2G(x) g'(x) + \frac{{g^{(4)}}(0)}{5} g^3 G > (ou < 0)$$ is equivalent to 
$$ \frac{g^2(x)}{2G(x)}+\frac{{g^{(4)}}(0)}{10}g^3(x) >(or <0) g'(x).$$ Moreover, in a neighborhood of $0$ one gets $$g^2(x)+\frac{{g^{(4)}}(0)}{5}g^3(x) G(x) = x^2+ \frac{{g^{(3)}}(0)}{3}x^4+...$$ 
$$\frac{2g G}{x}=x^2+\frac{{g^3}(0)}{12}x^4+...$$ Thus, according the hypothesis ${g^{(3)}}(0) <(or >) 0$ 
  $$ g^2(x)+\frac{{g^{(4)}}(0)}{5}g^3(x) G(x) <(or >) \frac{2g(x) G(x)}{x}$$  which is equivalent to $$\frac{g^2(x)}{2 G}+\frac{{g^{(4)}}(0)}{10}g^3(x)  <(or >) \frac{g(x) }{x}$$ which implying $ \frac{g(x)}{x} >(or <) g'(x)$ or equivalently $x (\frac{d}{dx}\frac{g(x)}{x})) < 0(or > 0)$.   

\bigskip

However, these results suppose generally the hypothesis $g^{(3)}(0)\neq 0$ holds.\\ Moreover, notice that neither [6] nor [2] have explicitly considered the case has $g^{(3)}(0)= 0$. Neverthless, we can deduce another consequence from Theorem 1-1. Indeed,
when $g''(0)=g^{(3)}(0)=0, g^{(4)}(0) < 0$, then condition 
$x (\frac{d}{dx}\frac{g}{x})) < 0(or > 0)$ falls 	and $ (\mathcal{C}_1): \ \frac{d}{dx}[\frac{G}{g^2}(x)] >(or <) \frac{{g^{(4)}}(0)}{5}G >(or <) \frac{d}{dx}[\frac{G}{g^2}(A(x))]$ appears to be the better monotonicity condition
for the period function $T(c)$ of (1).

\bigskip

\vspace{1cm}
\section {An example} Let us consider

$$g(x)=g_s(x)={\frac {1/2\, \left( x+s \right) \sinh \left( 2\,x \right) -1/4\,\cosh
 \left( 2\,x \right) +1/4}{s}}$$
$$g_s'(x) = {\frac { \left( x+s \right) \cosh \left( 2\,x \right) }{s}}$$
it is easy to see that $g_s$ verifies hypothesis $(\mathcal{H})$. The potential is then
$$G_s(x) = {\frac {1/4\,x\cosh \left( 2\,x \right) -1/4\,\sinh \left( 2\,x
 \right) +1/4\,s\cosh \left( 2\,x \right) +1/4\,x}{s}}-1/4
$$
The derivatives at $0$ are $$g''(0)=\frac{1}{s}, \ g^{4}(0)=\frac{12}{s}, \ g^{6}(0)=\frac{80}{s},...$$
Calculate the derivatives of the function $f(G)$ at $0$ one obtains 
$$f(0)= -\frac{1}{3s}, \ f'(0)= -\frac{7}{9}g''^3(0)+ \frac{g^{(4)}(0)}{5}=-{\frac {12}{5}}\,{s}^{-1}+{\frac {7}{9}}\,{s}^{-3},$$ $$ f''(0)=-\frac{28550(g''(0))^6-10320(g''(0))^3g^{4}(0)+81(g^{4}(0))^2+270(g''(0))g^{6}(0)}{4050 g''(0)}=$$ $${-\frac {1}{2025}}\,{\frac {16632\,{s}^{4}-14275-61920\,{s}^{2}}
{{s}^{5}}}$$

We have seen [3] that $$H(s,x)=g^2+(1/(3s)) g^3-2Gg'=g^2-2Gg'-f(0) g^3 $$ as a function of $x$ should change of sign for the value $s= .6447$. Thanks to {\it Maple} there is $x_0=-0.010737...$ such that \ $H(.647,x_0)= 0$. \\ 

Consider the following 
 $$H_1(s,x)=g^2+({\frac{1}{3s}})g^3-2Gg'-(\frac{-12}{5s}+\frac{7}{9s^3})g^3G$$ $$ H_1(s,x)=H(s,x)-f'(0)g^3 G= H(s,x)-(\frac{-12}{5s}+\frac{7}{9s^3})g^3 G$$
 
Here too $H_1(.647,x)$ should change of sign. Indeed, thanks to {\it Maple}   there is $x_1=-0.554537...$ such that \ $H_1(.647,x_1)= 0$. \\ 

Therefore, in order to prove the monotonicity of the period function we need to consider a better criteria. Let us consider the following

$$H_2(s,x)=g^2+(\frac{1}{3s})*g^3-2 G g'-[\frac{-12}{5s}+\frac{7}{9s^3}] g^3 G+\frac{[14275-61920 s^2+16632 s^4]}{4050 s^5} g^3 G^2$$
$$H_2(s,x)=H_1(s,x)-\frac{f''(0)}{2}g^3 G^2=H_1(s,x)+(\frac{[14275-61920 s^2+16632 s^4]}{4050 s^5} g^3 G^2$$

Thanks to {\it Maple} $H_2(s,x)$ should be negative for $x$ near $0$ and $s=.647$. That means condition $(\mathcal{C}_2)$ is satisfied while 
 $(\mathcal{C}_1)$ and $(\mathcal{C}_0)$ are not.\\ Thus, the energy-period function $T(c)$ is decreasing.

The detailled calculus are given below
 
\begin{figure}[th]
\begin{center}
\fbox{\includegraphics[width=\largeur,height=\hauteur]{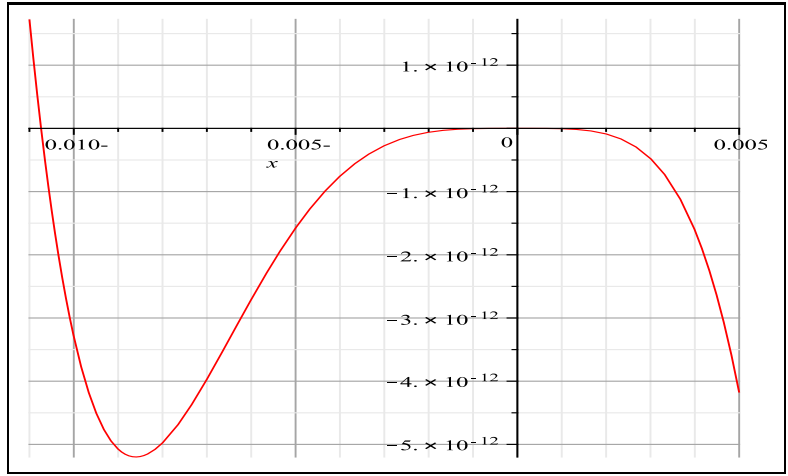}}
\end{center}
\caption{function $H(x)$ with $s=0.647,$ \ its zero is \ $ x_0=-0.010737... \cdots $\ This means  the
potential $G_{s}$ does not verify condition $(\mathcal{C}_{0})$ when $x \in [-\alpha, \alpha)]$ for $x_0 < \alpha$}
\label{fig1}%
\end{figure}

\begin{figure}[th]
\begin{center}
\fbox{\includegraphics[width=\largeur,height=\hauteur]{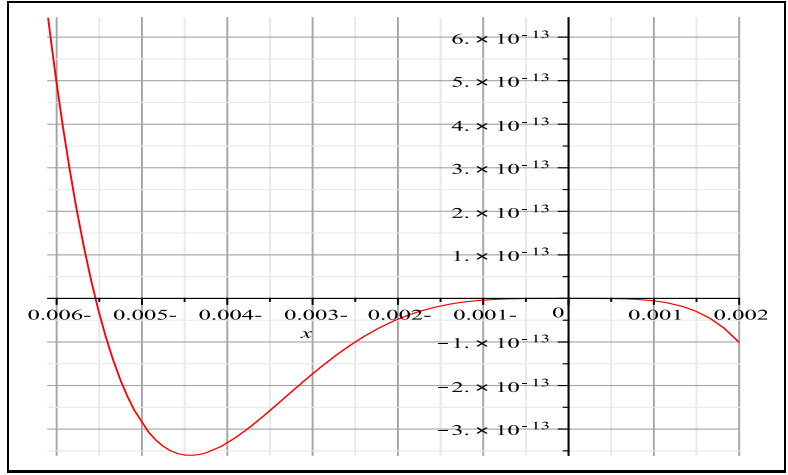}}
\end{center}
\caption{function $H_1(x)$ with $s=0.647,$ \ its zero is \ $ x_1=-0.005545... \cdots $ This means  the
potential $G_{s}$ does not verify condition $(\mathcal{C}_{1})$when $x \in [-\alpha, \alpha)]$ for $x_1 < \alpha$}
\label{fig2}%
\end{figure}

\begin{figure}[th]
\begin{center}
\fbox{\includegraphics[width=\largeur,height=\hauteur]{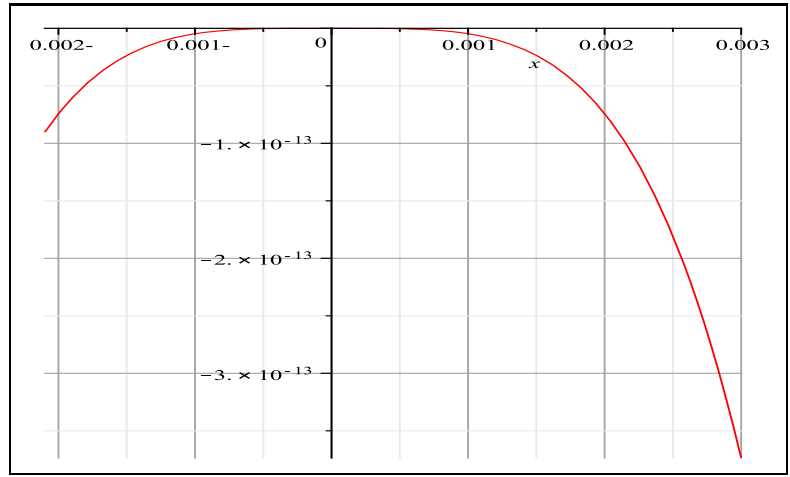}}
\end{center}
\caption{function $H_2(x)$ with $s=0.647.$ \ \ $H_2(x)<0$ for
$x \in [-\alpha,\alpha]$ and $x\neq0. $, here $\alpha = 0.04 \cdots$ This means $(\mathcal{C}_{2})$ can be satisfied}
\label{fig3}%
\end{figure}

\begin{figure}[th]
\begin{center}
\fbox{\includegraphics[width=\largeur,height=\hauteur]{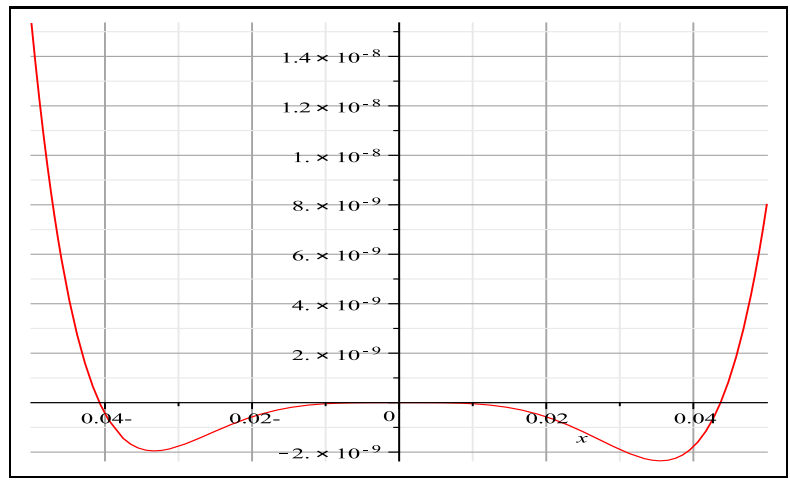}}
\end{center}
\caption{function $H_2(x)$ with $s=0.647.$ \ \ $H_2(x)$ has 2 zeros: \ \ \ \ \ $\alpha_1= -0.040743..$ 
and $\alpha_2=0.043699... $ This means  the
potential $G_{s}$ should verify condition $(\mathcal{C}_{2}),$ for $x \in [-\alpha,\alpha]$. So the period $T(c)$ is decreasing for $0<c<\bar c$}
\label{fig4}%
\end{figure}\clearpage

\section {Appendice: {\it Maple} computation details}

\bigskip

\begin{tabular}{|lll|}

$>$ & {\small eg:=(((x+s)/2)*sinh(2*x)-(cosh(2*x)-1)/4)/s : } \\ 
$>$ & {\small g:=unapply(eg,s,x);}\textsl{\ }\textit{\# function }$g_{s}$ \\ 
$>$ & {\small eG:=int(eg,x)-1/4:G:=unapply(eG,s,x);}\textrm{\ }\textit{\#
primitive of }$g_{s}$ \\ 
$>$ & {\small eg1:=diff(eg,x) : g1:=unapply(eg1,s,x);}\textrm{\ }\textit{\#
first derivated of }$g_{s}$ \\ 
$>$ & {\small eg2:=diff(eg1,x) : g2:=unapply(eg2,s,x);}\textrm{\ }\textit{\#
second derivated of }$g_{s}$ \\ 
$>$ & {\small eg3:=diff(eg2,x) : g3:=unapply(eg3,s,x);}\textrm{\ }\textit{\#
third derivated of }$g_{s}$ \\ 
$>$ & {\small eg4:=diff(eg3,x) : g4:=unapply(eg4,s,x);}\textrm{\ }\textit{\#
fourth derivated of }$g_{s}$ \\

$>$ & .\\
$>$ & .\\
$>$ & .\\

$>$ &  fsolve(H(.647, x), x = -0.15e-1 .. -0.1e-2)\\
                               -0.01073718589\\
$>$ &  fsolve(H1(.647, x), x = -0.15e-1 .. -0.1e-2)\\
                               -0.005545373709\\
$>$ &  fsolve(H2(.647, x), x = -0.5e-1 .. -0.4e-1)\\
                               -0.04074327315\\
$>$ &  fsolve(H2(.647, x), x = 0.4e-1 .. 0.5e-1)\\
                                0.04369965656\\
                                                            
 \end{tabular}                               
\vspace{1cm}

{\bf Coefficients of $g_s$}

\bigskip

Let us write $$g(x)=x+\frac{a_2}{2} x^2+\frac{a_3}{6} x^3+\frac{a_4}{24} x^4+\frac{a_5}{120} x^5+....$$ and $$G=\frac{1}{2} x^2+\frac{a_2}{6} x^3+\frac{a_3}{24} x^4+..., \  f(G)=b_0+b_1 G+\frac{b_2}{2} G^2+...$$

By identifying the coefficients we find after simplification
$${\it a_3}=5/3\,{{\it a_2}}^{2}, \quad {\it a_5}=-{\frac {140}{9}}\,{{\it a_2}}^{4}+7\,{\it a_2}\,{\it a_4}\quad 
{\it b_0}=-1/3\,{\it a_2}, $$ 
We thus obtain the first coefficients of the function $f$ 
$$ {\it b_1}=-1/5\,{\it a_4}+{\frac {7}{9}}\,{{\it a_2}}^{3}$$ 

$${\it b_2}=-\frac{28550(g''(0))^6-10320(g''(0))^3g^{4}(0)+81(g^{4}(0))^2+270(g''(0))g^{6}(0)}{2025 g''(0)}$$

\bigskip

{\bf Case of} $a_2=0$

This case is easier than the previous. We find after simplifying
$$ a_3 = a_5 = 0,\quad a_7 = (63/5)a_4^2,\quad a_9 = 66 a_6 a_4$$
Then it yields the first coefficients of $f$
 $$b_0=0,\quad b_1 = -(1/5)a_4,\quad b_2 = -(1/42)a_6,$$

$$  b_3 = -(1/810)a_8,\quad b_4 = -(1/27720)a_{10}+(13/600)a_4^3, $$ 

\newpage

{\bf References}

\bigskip

[1]\  A.R. Chouikha \quad {\it Period function and characterizations of Isochronous potentials}\quad arXiv:1109.4611, (2011).

\smallskip 
[2]\  A.R. Chouikha \quad {\it Monotonicity of the period function for some planar differential systems, I. Conservative and quadratic systems} \quad     Applic. Math., {\bf 32}, no. 3, p. 305-325, (2005). 

\smallskip  
[3]\  A.R. Chouikha and F. Cuvelier \quad {\it Remarks on some monotonicity conditions for the period function}\quad  Applic. Math., {\bf 26}, no. 3, p. 243-252, (1999).
                      
\smallskip 
[4]  \ S.N. Chow and D. Wang \quad {\it On the monotonicity of the period function of some second order equations}\quad Casopis Pest. Mat. {\bf 111}, p. 14-25, (1986).

\smallskip
[5]\ L.D. Landau E.M. Lifschitz \quad {\it Mechanics, Course of Theorical Physics} Vol 1, Pergamon Press, Oxford, (1960). 
                  
\smallskip 
[6] \ Z. Opial \quad {\it Sur les périodes des solutions de l'équation différentielle}\  $x'' + g(x) = 0$ \ Ann. Polon. Math., {\bf 10}, (1961), 49-72.

\smallskip
[7] \ F. Rothe \quad {\it Remarks on periods of planar Hamiltonian systems.} \quad SIAM J. Math. Anal., {\bf 24}, p.129-154, (1993).

\smallskip
[8] \ R. Schaaf \quad {\it A class of Hamiltonian systems with increasing periods} \quad J. Reine Angew. Math., {\bf 363}, p. 96-109, (1985). 

\end{document}